\documentclass[11pt]{amsart}
\usepackage{amssymb}
\usepackage{amsthm}
\usepackage{graphicx} 

\newtheorem{definition}{Definition}
\newtheorem{prop}{Proposition}

\newtheorem{thm}{Theorem}

\newtheorem{lemma}{Lemma}

\title[Every graph has an embedding with no non-hyperbolic knot]{Every graph has an embedding in $S^3$ containing no non-hyperbolic knot}
\author{Erica Flapan and Hugh Howards}

\begin{document}

\date \today
\maketitle

In contrast with knots, whose
properties depend only on their extrinsic topology in $S^3$, there is a rich 
interplay between the intrinsic structure of a graph and the 
extrinsic topology of all embeddings of the graph in $S^3$.   For example, it was shown in \cite{CG} that every embedding of the complete graph $K_7$ in $S^3$ contains a non-trivial knot.  Later in
 \cite{fl} it was shown that for every $m\in \mathbb N$, there is a complete graph $K_n$ such that every embedding of $K_n$ in $S^3$ contains a knot $Q$ (i.e., $Q$ is a subgraph of $K_n$) such that $|a_2(Q)|\geq m$, where $a_2$ is the second coefficient of the Conway polynomial of $Q$.  More recently, in \cite{fmn} it was shown that for every $m\in \mathbb N$, there is a complete graph $K_n$ such that every embedding of $K_n$ in $S^3$ contains a knot $Q$ whose minimal crossing number is at least $m$.  Thus there are arbitrarily complicated knots (as measured by $a_2$ and the minimal crossing number) in every embedding of a sufficiently large complete graph in $S^3$.  
  
 In light of these results, it is natural to ask whether there is a graph such that every embedding of that graph in $S^3$ contains a composite knot.  Or more generally, is there a graph such that every embedding of the graph in $S^3$ contains a satellite knot? Certainly, $K_7$ is not an example of such a graph since Conway and Gordon \cite{CG} exhibit an embedding of $K_7$ containing only the trefoil knot.   In this paper we answer this question in the negative.  In particular, we prove that every graph has an embedding in $S^3$ such that every non-trivial knot in that embedding is hyperbolic.  Our theorem implies that every graph has an embedding in $S^3$ which contains no composite or satellite knots.  By contrast, for any particular embedding of a graph we can add local knots within every edge to get an embedding such that every knot in that embedding is composite. 

Let $G$ be a graph.  There is an odd number $n$, such that $G$ is a minor of $K_n$.  We will show that for every odd number $n$, there is an embedding of $K_n$ in $S^3$ such that every non-trivial knot in that embedding of $K_n$ is hyperbolic.  It follows that there is an embedding of $G$ in $S^3$ which contains no non-trivial non-hyperbolic knots.

Let $n$ be a fixed odd number.  We begin by constructing a preliminary embedding of $K_n$ in $S^3$ as follows.  Let $h$ be a rotation of $S^3$ of order $n$ with fixed point set $\alpha\cong S^1$.  Let $V$ denote the complement of an open regular  neighborhood of the fixed point set $\alpha$.  Let $v_1$, \dots, $v_n$ be points in $V$ such that for each $i$, $h(v_i)=v_{i+1}$  (throughout the paper we shall consider our subscripts mod $n$).  These $v_i$ will be the vertices of the preliminary embedding of $K_n$.

\begin{definition}
By a {\bf solid annulus} we shall mean a 3-manifold with boundary which can be parametrized as $D\times I$ where $D$ is a disk. 
 We use the term {\bf the annulus boundary} of a solid annulus $D\times I$ to refer to the annulus $\partial D\times I$.  The {\bf ends} of $D\times I$ are the disks $D\times \{0\}$ and $D\times \{1\}$.  If $A$ is an arc in a solid annulus $W$ with one endpoint in each end of $W$, and $A$ co-bounds a disk in $W$ together with an arc in $\partial W$, then we say that  $A$ is a {\bf longitudinal arc} of $W$.

\end{definition}

\begin{figure}[h]
\includegraphics{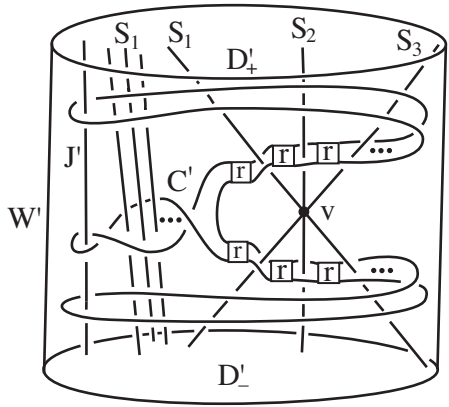}
\caption {For each $i$, $W'\cap S_i$ is a collection of $n+i$ untangled longitudinal arcs.}
\label{twists}
\end{figure}

  As follows, we embed the edges of $K_n$ as simple closed curves in the quotient space $S^3/h=S^3$. Observe that since $V$ is a solid torus, $V'=V/h$ is also a solid torus.  Let $D'$ denote  a meridional disk for $V'$ which does not contain the point $v=v_1/h$.  Let $W'$ denote the solid annulus $\mathrm{cl}(V'-D')$ with ends $D'_+$ and $D'_-$.   Since $n$ is odd, we can choose unknotted simple closed curves $S_1$, \dots, $S_{\frac{n-1}{2}}$ in the solid torus $V'$ such that each $S_i$ contains $v$ and has winding number $n+i$ in $V'$, the $S_i$ are pairwise disjoint except at $v$, and for each $i$, $W'\cap S_i$  is a collection of $n+i$ untangled longitudinal arcs (see Figure \ref{twists}).    
  
We define two additional simple closed curves $J'$ and $C'$ in $V'$  whose intersections with $W'$ are illustrated in Figure \ref{twists} as follows.  First, choose a simple closed curve $J'$ in $V'$, whose intersection with $W'$ is a longitudinal arc which is disjoint from and untangled with $S_1\cup\dots \cup S_{\frac{n-1}{2}}$.  Next we let $C'$ be the unknotted simple closed curve in 
 Ê$W'-(S_1\cup\dots \cup S_{\frac{n-1}{2}} \cup J')$ whose projection is illustrated in Figure \ref{twists}. ÊIn particular, $CÕ$ contains one half 
 twist between $J'$ and the set of arcs of $S_1\cup\dots \cup 
 S_{\frac{n-1}{2}}$ which do not contain $v$, another half twist between 
 those arcs of $S_1\cup\dots \cup S_{\frac{n-1}{2}}$ and the set of arcs 
 containing $v$, and $r$ full-twists between each of the individual 
 Êarcs of Ê$S_i$ and $S_{i+1}$ containing $v$. ÊWe will determine the value of $r$ later.

Each of the $\frac{n-1}{2}$ simple closed curves $S_1$, \dots, $S_\frac{n-1}{2}$ lifts to a simple closed curve consisting of $n$ consecutive edges of $K_n$.  The vertices $v_1$, \dots, $v_n$ together with these $\frac{n(n-1)}{2}$ edges gives us a preliminary embedding $\Gamma_1$ of $K_n$ in $S^3$.  

Lift the meridional disk $D'$ of the solid torus $V'$ to $n$ disjoint meridional disks $D_1$, \dots, $D_n$ of the solid torus $V$.  Lift the simple closed curve $C'$ to $n$ disjoint simple closed curves $C_1$, \dots, $C_n$, and lift the simple closed curve $J'$ to $n$ consecutive arcs $J_1$, \dots, $J_n$ whose union is a simple closed curve $J$.  The closures of the components of $V-(D_1\cup \dots \cup D_n)$ are solid annuli, which we denote by $W_1$, \dots, $W_n$.  The subscripts of all of the lifts are chosen consistently so that for each $i$, $v_i\in W_i$, $C_i\cup J_i\subseteq W_i$, and $D_i$ and $D_{i+1}$ are the ends of the solid annulus $W_i$. For each $i$, the pair $(W_i-(C_i \cup J_i),(W_i-(C_i \cup J_i))\cap \Gamma_1))$ is homeomorphic to $(W'-(C' \cup J'),(W'-(C' \cup J'))\cap(S_1\cup\dots \cup S_{\frac{n-1}{2}}))$.  For each $i$, the solid annulus $W'$ contains $n+i-1$ arcs of $S_i$ which are disjoint from $v$.  Hence each edge of the embedded graph $\Gamma_1$ meets each solid annulus $W_i$ in at least one arc not containing $v_i$.

\begin{figure}[h]
\includegraphics{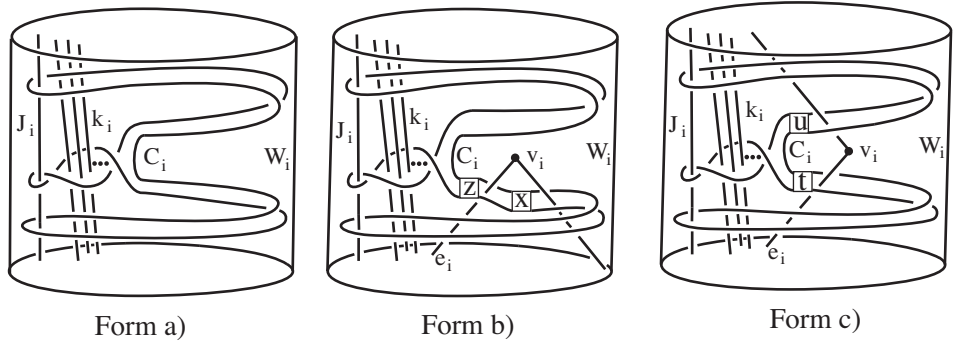}
\caption{The forms of  $(W_i,C_i \cup J_i\cup k_i\cup e_i)$.}
\label{WicapK}
\end{figure}

   Let $\kappa$ be a simple closed curve in $\Gamma_1$.   For each $i$, we let $k_i$ denote the set of those arcs of $\kappa\cap W_i$ which do not contain $v_i$, and let $e_i$ denote either the single arc of $\kappa\cap W_i$ which does contain $v_i$ or the empty set if $v_i$ is not on $\kappa$.  Observe that since $\kappa$ is a simple closed curve, it contains at least three edges of $\Gamma_1$; and as we observed above, each edge of $\kappa$ contains at least one arc of $k_i$.   Thus for each $i$, $k_i$ contains at least three arcs.   Either $e_i$ is empty, the endpoints of $e_i$ are in the same end of the solid annulus $ W_i$, or the endpoints of $e_i$ are in different ends of $W_i$.  We illustrate these three possibilities for  $(W_i,C_i \cup J_i\cup k_i\cup e_i)$ In Figure \ref{WicapK} as forms a), b) and c) respectively.  The number of full-twists represented by the labels $t$, $u$, $x$, or $z$ in Figure \ref{WicapK} is some multiple of $r$ depending on the particular simple closed curve $\kappa$. 

\begin{figure}[h]
\includegraphics{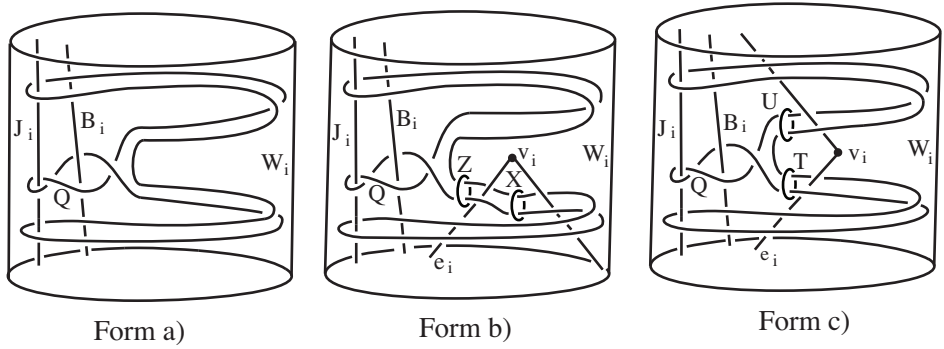}
\caption{The forms of $W_i$ with associated simple closed curves and the arc $B_i$.}
\label{curves}
\end{figure}

For each of the forms of $(W_i,C_i \cup J_i\cup k_i\cup e_i)$ illustrated in Figure \ref{WicapK}, we will associate an additional arc and an additional collection of simple closed curves as follows (illustrated in Figure \ref{curves}).   Let the arc $B_i$ be the core of a solid annulus neighborhood of the union of the arcs $k_i$ in $W_i$ such that $B_i$ is disjoint from $J_i$, $C_i$, and $e_i$.  Let the simple closed curve $Q$ be obtained from $C_i$ by removing the full twists $z$, $x$, $t$, and $u$.  Let $Z$, $X$, $T$, and $U$ be unknotted simple closed curves which wrap around $Q$ in place of $z$, $x$, $t$, and $u$ as illustrated in Figure \ref{curves}.

ÊÊFor each $i$, let $M_i$ denote an unknotted solid torus in $S^3$ 
 Êobtained by gluing together two identical copies of $W_i$ along $D_i$ and $D_{i+1}$, 
 Êmaking sure that the end points of the arcs of $J_i$, $B_i$, and $e_i$ match 
 Êup with their counterparts in the second copy to get simple closed 
 curves $j$, $b$, and $E$ respectively in $M_i$. ÊThus $M_i$ has a $180^{\circ}$ 
 Êrotational symmetry around a horizontal line which goes through the center of the figure and 
 Êthe end points of both copies of $J_i$, $B_i$, and $e_i$.   Recall that in form a), $e_i$ is the empty set, and hence so is $E$.  Let $Q_1$ and $Q_2$, $X_1$ and $X_2$, $Z_1$ and $Z_2$, $T_1$ and $T_2$, and $U_1$ and $U_2$ denote the doubles of the unknotted simple closed curves $Q$, $X$, $Z$, $T$, $U$ respectively.  
  
\begin{figure}[h]
\includegraphics{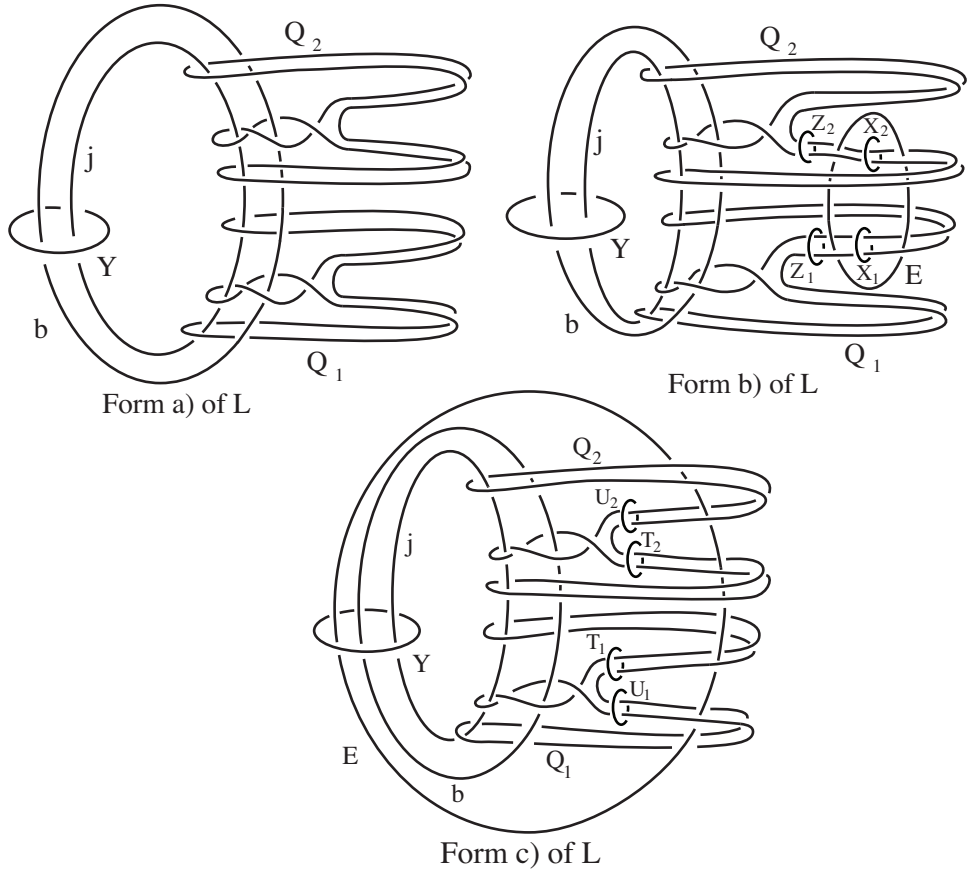}
\caption{The possible forms of  the link $L$.}
\label{link}
\end{figure}

    Let $Y$ denote the core of the solid torus $\mathrm{cl}(S^3-M_i)$. We associate to Form a) of Figure \ref{curves} the link $L=Q_1\cup Q_2\cup j\cup b\cup Y$.  We associate to Form b) of Figure \ref{curves} the link $L=   Q_1\cup Q_2\cup j\cup b\cup Y\cup E\cup X_1\cup X_2\cup Z_1\cup Z_2$.  We associate to Form c) of Figure \ref{curves} the link $L= Q_1\cup Q_2\cup j\cup b\cup Y\cup E\cup T_1\cup T_2\cup U_1\cup U_2$.    Figure \ref{link} illustrates the three forms of the link
$L$.

  The software program SnapPea (available at http://www.geometrygames.\linebreak
  org/SnapPea/index.html) can be used to determine whether or not a given knot or link in $S^3$ is hyperbolic, and if so SnapPea estimates the hyperbolic volume of the complement.  We used SnapPea to verify that each of the three forms of the link $L$ illustrated in Figure \ref{link} is hyperbolic.  
  
A 3-manifold is unchanged by doing Dehn surgery on an unknot if the boundary slope of the surgery is the reciprocal of an integer (though such surgery may change a knot or link in the manifold).   According to Thurston's Hyperbolic Dehn Surgery Theorem \cite{BP,Th}, all but finitely many Dehn fillings of a hyperbolic link complement yield a hyperbolic manifold.   Thus there is some $r\in \mathbb{N}$ such that for any $m\geq r$,  if we do Dehn filling with slope $\frac{1}{m}$ along the components $X_1$, $X_2$, $Z_1$, $Z_2$ of the link $L$ in form b) or along the components $T_1$, $T_2$, $U_1$, $U_2$ of the link $L$ in form c), then we obtain a hyperbolic link $\overline{Q}_1\cup \overline{Q}_2\cup j\cup b\cup Y \cup E$, where the simple closed curves $\overline{Q}_1$ and $\overline{Q_2}$ are obtained by adding $m$ full twists to $Q_1$ and $Q_2$ in place of each of the surgered curves.

We fix the value of $r$ according to the above paragraph, and this is the value of $r$ that we use in Figure \ref{twists}.  Recall that the number of twists $x$, $z$, $u$, and $t$ in the simple closed curves $C_i$ in Figure \ref{WicapK} are each a multiple of $r$. Thus the particular simple closed curves $C_i$ are determined by our choice of $r$ together with our choice of the simple closed curve $\kappa$.   Now we do Dehn fillings along $X_1$ and $X_2$ with slope $\frac{1}{x}$, along $Z_1$ and $Z_2$ with slope $\frac{1}{z}$, along $U_1$ and $U_2$ with slope $\frac{1}{u}$, and along $T_1$ and $T_2$ with slope $\frac{1}{t}$.  Since $x$, $z$, $u$, and $t$ are each greater than or equal to $r$, the link $\overline{Q}_1\cup \overline{Q}_2\cup j\cup b\cup Y \cup E$ that we obtain will be hyperbolic.  In form a), $E$ is the empty set and the link $Q_1\cup Q_2\cup j\cup b\cup Y \cup E$ was already seen to be hyperbolic using SnapPea.  In this case, we do no surgery and we let the simple closed curves $\overline{Q}_1=Q_1$ and $\overline{Q}_2=Q_2$.  It follows that each form of $M_i-(\overline{Q}_1\cup \overline{Q}_2\cup j\cup b\cup E)$ is a hyperbolic 3-manifold.   Observe that $M_i-(\overline{Q}_1\cup \overline{Q}_2\cup j\cup b\cup E)$ is the double of $W_i-(C_i\cup J_i\cup B_i\cup e_i)$.
 
  Now that we have fixed $C_i$, we let $N(C_i)$, $N(J_i)$, $N(B_i)$, and $N(e_i)$ be pairwise disjoint regular neighborhoods of $C_i$, $J_i$, $B_i$, and $e_i$ respectively in the interior of each of the forms of the solid annulus $W_i$ (illustrated in Figure \ref{WicapK}).  We choose $N(B_i)$ such that it contains the union of the arcs $k_i$.  Note that in form a) $e_i$ is the empty set and hence so is $N(e_i)$.  Let $N(k_i)$ denote a collection of pairwise disjoint regular neighborhoods one containing each arc of $k_i$ such that $N(k_i)\subseteq N(B_i)$.  Let $V_i=\mathrm{cl}(W_i-(N(C_i) \cup N(J_i)\cup N(B_i)\cup N(e_i)))$, let $\Delta=\mathrm{cl}(N(B_i)-N(k_i))$, and let $V_i'=V_i\cup\Delta$.  Since $N(B_i)$ is a solid annulus, it has a product structure $D^2\times I$.  Without loss of generality, we assume that each of the components of $N(k_i)$ respects the product structure of $N(B_i)$.  Thus $\Delta=F\times I$ where $F$ is a disk with holes.  

\begin{definition}  Let  $X$ be a 3-manifold.  A sphere in $X$ is said to be {\bf essential} if it does not bound a ball in $X$. A properly embedded disk $D$ in $X$ is said to be {\bf essential} if $\partial D$ does not bound a disk in $\partial X$. A properly embedded annulus is said to be {\bf essential} if it is incompressible and not boundary parallel.
A torus in $X$ is said to be {\bf essential} if it is incompressible and not boundary parallel.  
\end{definition}
  
\begin{lemma}\label{L:annuli}
For each $i$, $V_i'$ contains no essential torus, sphere, or disk whose boundary is in $D_i\cup D_{i+1}$.  Also, any incompressible annulus in $V_i'$ whose boundary is in $D_i\cup D_{i+1}$ is either boundary parallel or can be expressed as $\sigma\times I$ (possibly after a change in parameterization of $\Delta$), where $\sigma$ is a non-trivial simple closed curve in $D_i\cap\Delta$. 
\end{lemma}

\begin{proof}    Since $k_i$ contains at least three disjoint arcs,  $F$ is a disk with at least three holes.  Let $\beta$ denote the double of $\Delta$ along $\Delta\cap(D_i\cup D_{i+1})$.  Then $\beta=F\times S^1$.  Now it follows from Waldhausen \cite{Wa} that $\beta$ contains no essential sphere or properly embedded disk, and any incompressible torus in $\beta$ can be expressed as $\sigma\times S^1$ (after a possible change in parameterization of $\beta$) where $\sigma$ is a non-trivial simple closed curve in $D_i\cap\Delta$.

Let $\nu$ denote the double of $V_i$ along $V_i\cap(D_i\cup D_{i+1})$.  Observe that $\nu\cup \beta$ is the double of $V_i'$ along $V_i'\cap (D_i\cup D_{i+1})$.  Now the interior of $\nu$ is homeomorphic to $M_i-(\overline{Q}_1\cup \overline{Q}_2\cup j\cup b\cup E)$.  Since we saw above that $M_i-(\overline{Q}_1\cup \overline{Q}_2\cup j\cup b\cup E)$ is hyperbolic, it follows from Thurston \cite{Th,Th2} that $\nu$ contains no essential sphere or torus, or properly embedded disk or annulus.  

We see as follows that $\nu\cup\beta$ contains no essential sphere and any essential torus in $\nu\cup\beta$ can be expressed (after a possible change in parameterization of $\beta$) as $\sigma\times S^1$, where $\sigma$ is a non-trivial simple closed curve in $D_i\cap\Delta$.  Let $\tau$ be an essential sphere or torus in $\nu\cup\beta$, and let $\gamma$ denote the torus $\nu\cap \beta$. By doing an isotopy as necessary, we can assume that $\tau$ intersects $\gamma$ in a minimal number of disjoint simple closed curves.  Suppose there is a curve of intersection which bounds a disk in the essential surface $\tau$.  Let $c$ be an innermost curve of intersection on $\tau$ which bounds a disk $\delta$ in $\tau$.   Then $\delta$ is a properly embedded disk in either $\gamma$ or $\beta$.  Since neither $\nu$ nor $\beta$ contains a properly embedded essential disk or an essential sphere, there is an isotopy of $\tau$ which removes $c$ from the collection of curves of intersection.  Thus by the minimality of the number of curves in $\tau\cap\gamma$, we can assume that none of the curves in $\tau\cap\gamma$ bounds a disk in $\tau$.  

Suppose that $\tau$ is an essential sphere in $\nu\cup\beta$.  Since none of the curves in $\tau\cap\gamma$ bounds a disk in $\tau$, $\tau$ must be contained entirely in either $\nu$ or $\beta$.  However, we saw above that neither $\nu$ nor $\beta$ contains any essential sphere. Thus $\tau$ cannot be an essential sphere, and hence must be an essential torus.  Since $\tau\cap\gamma$ is minimal,  if $\tau\cap\nu$ is non-empty, then the components of $\tau$ in $\nu$ are all incompressible annuli.  However, we saw above that $\nu$ contains no essential annuli.   Thus $\tau\cap\nu$ is empty.  Since $\nu$ contains no essential torus, the essential tori $\tau$ must be contained in $\beta$.  Hence $\tau$ can be expressed (after a possible change in parameterization of $\beta$) as $\sigma\times S^1$, where $\sigma$ is a non-trivial simple closed curve in $D_i\cap\Delta$. 

Now we consider essential surfaces in $V_i'$.  Suppose that $V_i'$ contains an essential sphere $S$.  Since $\nu\cap\beta$ contains no essential sphere, $S$ bounds a ball $B$ in $\nu\cap \beta$.  Now the ball $B$ cannot contain any of the boundary components of $\nu\cap \beta$.  Thus $B$ cannot contain either $D_i$ or $D_{i+1}$.  Since $S$ is disjoint from $D_i\cup D_{i+1}$, it follows that $B$ must be disjoint from $D_i\cup D_{i+1}$.  Thus $B$ is contained in $V_i'$.  Hence $V_i'$ cannot contain an essential sphere.

We see as follows that $V_i'$ cannot contain an essential disk whose boundary is in $D_i\cup D_{i+1}$.  Let $\epsilon$ be a disk in $V_i'$ whose boundary is in $D_i\cup D_{i+1}$.   Let $\epsilon'$ denote the double of $\epsilon$ in $\nu\cup\beta$.  Then $\epsilon'$ is a sphere which meets $D_i\cup D_{i+1}$ in the simple closed curve $\partial \epsilon$.  Since $\nu\cup \beta$ contains no essential sphere, $\epsilon'$ bounds a ball $B$ in $\nu\cup \beta$.  It follows that $B$ cannot contain any of the boundary components of $\nu\cup \beta$.  Thus $B$ cannot contain any of the boundary components of $D_i\cup D_{i+1}$.  Therefore, $D_i\cup D_{i+1}$ intersects the ball $B$ in a disk bounded by $\partial\epsilon$.  Hence the simple closed curve $\partial\epsilon$ bounds a disk in $(D_i\cup D_{i+1})\cap V_i'$, and therefore the disk $\epsilon$ was not essential in $V_i'$.  Thus,  $V_i'$ contains no essential disk whose boundary is in $D_i\cup D_{i+1}$.

 Now suppose that $V_i'$ contains an essential torus $T$.  Suppose that $T$ is not essential in $\nu\cup \beta$.  Then either $T$ is boundary parallel or $T$ is compressible in $\nu\cup\beta$.  However, $T$ cannot be boundary parallel in $\nu\cup\beta$ since $T\subseteq V_i'$.  Thus $T$ must be compressible in $\nu\cup\beta$.  Let $\delta$ be a compression disk for $T$ in $\nu\cup\beta$.  Since $V_i'$ contains no essential sphere or essential disk whose boundary is in $D_i\cup D_{i+1}$, we can use an innermost disk argument to push $\delta$ off of $D_i\cup D_{i+1}$.  Hence $T$ is compressible in $V_i'$, contrary to our initial assumption.  Thus $T$ must be essential in $\nu\cup \beta$.  It follows that $T$ has the form $\sigma\times S^1$, where $\sigma\subseteq D_i\cap \Delta$.  However, since $\nu\cup \beta$ is the double of $V_i'$, the intersection of $\sigma\times S^1$ with $V_i'$ is an annulus $\sigma\times I$.  In particular, $V_i'$ cannot contain  $\sigma\times S^1$.  Therefore, $V_i'$ cannot contain an essential torus.

Suppose that $V_i'$ contains an incompressible annulus $\alpha$ whose boundary is in $D_i\cup D_{i+1}$.  Let $\tau$ denote the double of $\alpha$ in $\nu\cup\beta$.  Then $\tau$ is a torus.  If $\tau$ is essential in $\nu\cup\beta$, then we saw above that $\tau$ can be expressed as $\sigma\times S^1$ (after a possible change in parameterization of $\beta$) where $\sigma$ is a non-trivial simple closed curve in $D_i\cap\Delta$.  In this case, $\alpha$ can be expressed as $\sigma\times I$. 

On the other hand, if $\tau$ is inessential in $\nu\cup \beta$, then either $\tau$ is parallel to a component of $\partial (\nu\cup\beta)$, or $\tau$ is compressible in $\nu\cup\beta$.   If $\tau$ is parallel to a boundary component of $\nu\cup\beta$, then $\alpha$ is parallel to the annulus boundary component of $W_i$, $N(J_i)$, $N(e_i)$, $N(B_i)$, or one of the boundary components of $N(k_i)$. 

Thus we suppose that the torus $\tau$ is compressible in $\nu\cup\beta$.  In this case, it follows from an innermost loop outermost arc argument that either the annulus $\alpha$ is compressible in $V_i'$ or $\alpha$ is $\partial$-compressible in $V_i'$.  Since we assumed $\alpha$ was incompressible in $V_i'$, $\alpha$ must be $\partial$-compressible in $V_i'$.   Now according to a lemma of Waldhausen \cite{Wa}, if a 3-manifold contains no essential sphere or properly embedded essential disk, then any annulus which is incompressible but boundary compressible must be boundary parallel.  We saw above that $V_i'$ contains no essential sphere or essential disk whose boundary is in $D_i\cup D_{i+1}$.  Since the boundary of the incompressible annulus $\alpha$ is contained in $D_i\cup D_{i+1}$, it follows from Waldhausen's Lemma that $\alpha$ is boundary parallel in $V_i'$.   \end{proof}  
    
\bigskip

It follows from Lemma \ref{L:annuli} that for any $i$, any incompressible annulus in $V_i'$ whose boundary is in $D_i\cup D_{i+1}$ is either parallel to an annulus in $D_i$ or $D_{i+1}$ or co-bounds a solid annulus in the solid annulus $W_i$ with ends in $D_i\cup D_{i+1}$.  Recall that $\kappa$ is a simple closed curve in $\Gamma_1$ such that $\kappa\cap W_i= k_i\cup e_i$.  Also $J=J_1\cup\dots\cup J_n$.  Let $N(\kappa)$ and $N(J)$ be regular neighborhoods of the simple closed curves $\kappa$ and $J$ respectively, such that for each $i$, $N(\kappa)\cap W_i=N(k_i)\cup N(e_i)$, and $N(J)\cap W_i=N(J_i)$.  Recall that $V=W_1\cup \dots \cup W_n$.  Thus $\mathrm{cl}(V-(N(C_1)\cup \dots\cup N(C_n) \cup N(J) \cup N(\kappa))=V_1'\cup\dots\cup V_n'$.
\medskip

\begin{prop}\label{prop:notori}    $H=\mathrm{cl}(V-(N(C_1)\cup \dots\cup N(C_n) \cup N(J) \cup N(\kappa))$ contains no essential sphere or torus.  
\end{prop}

\begin{proof}  Suppose that $S$ is an essential sphere in $H$.  Without loss of generality, $S$ intersects the collection of disks $D_1\cup\dots\cup D_n$ transversely in a minimal number of simple closed curves. By Lemma \ref{L:annuli}, for each $i$,  $V_i'$ contains no essential sphere or essential disk whose boundary is in $D_i\cup D_{i+1}$.  Thus the sphere $S$ cannot be entirely contained in one $V_i'$.  Let $c$ be an innermost curve of intersection on $S$.  Then $c$ bounds a disk $\delta$ in some $V_i'$.  However, since the number of curves of intersection is minimal, $\delta$ must be essential, contrary to Lemma \ref{L:annuli}.  Hence $H$ contains no essential sphere.

 Suppose $T$ is an incompressible torus in $H$.  We show as follows that $T$ is parallel to some boundary component of $H$. Without loss of generality, the torus $T$ intersects the collection of disks $D_1\cup\dots\cup D_n$ transversely in a minimal number of simple closed curves. By Lemma \ref{L:annuli}, for each $i$,  $V_i'$ contains no essential torus, essential sphere, or essential disk whose boundary is in $D_i\cup D_{i+1}$.  Thus the torus $T$ cannot be entirely contained in one $V_i'$.  Also, by the minimality of the number of curves of intersection, we can assume that if $V_i'\cap T$ is nonempty, then it consists of a collection of incompressible annuli in $V_i'$ whose boundary components are in $D_i\cup D_{i+1}$.  Furthermore,  by Lemma \ref{L:annuli}, each such annulus is either boundary parallel or is contained in $N(B_i)$ and can be expressed (after a possible change in parameterization of $N(B_i)$) as $\sigma_i\times I$ for some non-trivial simple closed curve  $\sigma_i$ in $D_i\cap \Delta$.  If some annulus component of $V_i'\cap T$ is parallel to an annulus in $D_i\cup D_{i+1}$, then we could remove that component by an isotopy of $T$.  Thus we can assume that each annulus in $V_i'\cap T$ is parallel to the annulus boundary component of one of the solid annuli $W_i$, $N(J_i)$, or $N(e_i)$, or can be expressed as $\sigma_i\times I$.  In any of these cases the annulus co-bounds a solid annulus in $W_i$ with ends in $D_i\cup D_{i+1}$.

Consider  some $i$, such that $V_i'\cap T$ is non-empty.  Hence it contains an incompressible annulus $A_i$ which has one of the above forms.  By the connectivity of the torus $T$, either there is an incompressible annulus $A_{i+1}\subseteq V_{i+1}'\cap T$ such that $A_i$ and $A_{i+1}$ share a boundary component, or there is an incompressible annulus $A_{i-1}\subseteq V_{i-1}'\cap T$, such that $A_i$ and $A_{i-1}$ share a boundary component, or both.  We will assume, without loss of generality, that there is an incompressible annulus $A_{i+1}\subseteq V_{i+1}'\cap T$ such that $A_i$ and $A_{i+1}$ share a boundary component. Now it follows that $A_i$ co-bounds a solid annulus $F_i$ in $W_i$ with ends in $D_i\cup D_{i+1}$, and $A_{i+1}$ co-bounds a solid annulus $F_{i+1}$ in $W_{i+1}$ together with two disks in $D_{i+1}\cup D_{i+2}$.  Hence the solid annuli $F_i$ and $F_{i+1}$ meet in one or two disks in $D_{i+1}$.  

We consider several cases where $A_i$ is parallel to some boundary component of $V_i'$.  Suppose that $A_i$ is parallel to the annulus boundary component of the solid annulus $N(J_i)$.   Then the solid annulus $F_i$ contains $N(J_i)$ and is disjoint from the arcs $k_i$ and $e_i$.  Now the arcs $J_i$ and $J_{i+1}$ share an endpoint contained in $F_i\cap F_{i+1}$, and there is no endpoint of any arc of $k_i$ or $e_i$ in $F_i\cap F_{i+1}$.  It follows that the solid annulus $F_{i+1}$ contains the arc $J_{i+1}$ and contains no arcs of $k_{i+1}$.   Hence by Lemma \ref{L:annuli}, the incompressible annulus $A_{i+1}$ must be parallel to $\partial N(J_{i+1})$.  Continuing from one $V_i'$ to the next, we see that in this case, $T$ is parallel to $\partial N(J)$.  

Suppose that $A_i$ is parallel to the annulus boundary component of the solid annulus $\partial N(e_i)$ or one of the solid annuli in $\partial N(k_i)$.  Using an argument similar to the above paragraph, we see that $A_{i+1}$ is parallel to the annulus boundary component of the solid annulus $\partial N(e_{i+1})$ or one of the solid annuli in $\partial N(k_{i+1})$.  Continuing as above, we see that in this case $T$ is parallel to $\partial N(\kappa)$.

Suppose that the annulus $A_i$ is parallel to the annulus boundary component of the solid annulus $W_i$.  Then the solid annulus $F_i$ contains all of the arcs of $J_i$, $k_i$, and $e_i$.  It follows as above that the solid annulus $F_{i+1}$ contains the arc $J_{i+1}$ and some arcs of $k_{i+1}\cup e_{i+1}$.  Thus by Lemma \ref{L:annuli}, $A_{i+1}$ must be parallel to the annulus boundary component of the solid annulus $W_{i+1}$.  Continuing in this way, we see that in this case $T$ is parallel to $\partial V$.

Thus we now assume that no component of any $V_i'\cap T$ is parallel to an annulus boundary component of  $V_i'$.  Hence if any $V_i'\cap T$ is non-empty, then by Lemma \ref{L:annuli}, $V_i'\cap T$ consists of disjoint incompressible annuli in $N(B_i)$ which can each be expressed (after a possible re-parametrization of $N(B_i)$) as $\sigma_i\times I$ for some non-trivial simple closed curve  $\sigma_i\subseteq D_i\cap \Delta$.   Choose $i$ such that $V_i'\cap T$ is non-empty.  Since $N(B_i)$ is a solid annulus, there is an innermost incompressible annulus $A_i$ of $N(B_i)\cap T$.  Now $A_i$ bounds a solid annulus $F_i$ in $N(B_i)$, and $F_i$ contains more than one arc of $k_i$. Since $A_i$ is innermost in $N(B_i)$, $\mathrm{int}(F_{i})$ is disjoint from $T$.  Now there is an incompressible annulus $A_{i+1}$ in $V_{i+1}'\cap T$, such that $A_i$ and $A_{i+1}$ meet in a circle in $D_{i+1}$.  Furthermore, this circle bounds a disk in $D_{i+1}$ which is disjoint from $T$, and by our assumption is contained in $N(B_i)$.  Thus by Lemma  \ref{L:annuli}, the incompressible annulus $A_{i+1}$ has the form $\sigma_{i+1}\times I$ for some non-trivial simple closed curve  $\sigma_{i+1}\subseteq D_{i+1}\cap \Delta$  .  Thus $A_{i+1}$ bounds a solid annulus $F_{i+1}$ in $N(B_{i+1})$, and $\mathrm{int}(F_{i+1})$ is also disjoint from $T$.  We continue in this way considering consecutive annuli to conclude that for every $j$, every component $A_j$ of $T\cap V_j'$ is an incompressible annulus which bounds a solid annulus $F_j$ whose interior is disjoint from $T$.  

Recall that $V=W_1\cup \dots \cup W_n$ is a solid torus.  Let $Q$ denote the component of $V-T$ which is disjoint from $\partial V$.  Then $Q$ is the union of the solid annuli $F_j$. Since some $F_i$ contains some  arcs of $k_i$, the simple closed curve $\kappa$ must be contained in $Q$. 

 Recall that the simple closed curve $\kappa$ contains at least three vertices of the embedded graph $\Gamma_1$.  Also each vertex of $\kappa$ is contained in some arc $e_j$.  Since each such $e_j\subseteq \kappa\subseteq Q$,  some component $F_j$ of $Q\cap W_j$ contains the arc $e_j$.  By our assumption, for any $V_i'\cap T$ which is non-empty, $V_i'\cap T$ consists of disjoint incompressible annuli in $N(B_i)$.  In particular, $V_j\cap T\subseteq N(B_i)$.  Now the annulus boundary of $F_j$ is contained in $N(B_j)$, and hence $F_j\subseteq N(B_j)$.   But this is impossible since $e_j\subseteq F_j$ and $e_j$ is disjoint from $N(B_j)$.  Hence our assumption that no component of any $V_i'\cap T$ is parallel to an annulus boundary component of  $V_i'$ is wrong.  Thus, as we saw in the previous cases, $T$ must be parallel to a boundary component of $H$.  Therefore $H$ contains no essential annulus.  \end{proof}
\bigskip

Recall that the value of $r$, the simple closed curves, and the manifold $H$, all depend on the particular choice of simple closed curve $\kappa$.  In the following theorem, we do not fix a particular $\kappa$, so none of the above are fixed.

\begin{thm}  
Every graph can be embedded in $S^3$ in such a way that every non-trivial knot in the embedded graph is hyperbolic.
\end{thm}

\begin{proof}  Let $G$ be a graph, and let $n\geq 3$ be an odd number such that $G$ is a minor of the complete graph  on $n$ vertices $K_n$. Let $\Gamma_1$ be the embedding of $K_n$ given in our preliminary construction.  Then, $\Gamma_1$ contains at most finitely many simple closed curves, $\kappa_1$, \dots, $\kappa_m$.  For each $\kappa_j$, we use Thurston's Hyperbolic Dehn Surgery Theorem \cite{BP, Th} to choose an $r_j$ in the same manner that we chose $r$ after we fixed a particular simple closed curve $\kappa$.  Now let $R=\mathrm{max}\{r_1,\dots,r_m\}$, and let $R$ be the value of $r$ in Figure \ref{twists}.  This determines the simple closed curves $C_1$, \dots, $C_n$.

Let $P=\mathrm{cl}(V-(N(C_1)\cup \dots\cup N(C_n) \cup N(J )))$ where $V$ and $J$ are given in our preliminary construction.  Then the embedded graph $\Gamma_1\subseteq P$.  For each $j=1$, \dots $m$, let $H_j=\mathrm{cl}(P-N(\kappa_j))$.  It follows from Proposition \ref{prop:notori} that each $H_j$ contains no essential sphere or torus. Since each $H_j$ has more than three boundary components, no $H_j$ can be Seifert fibered.  Hence by Thurston's Hyperbolization Theorem \cite{Th2},  every $H_j$ is a hyperbolic manifold.

 We will glue solid tori $Y_1$, \dots, $Y_{n+2}$ to $P$ along its $n+2$ boundary components $\partial V$, $\partial N(C_1)$, \dots, $\partial N(C_n)$, and $\partial N(J)$ to obtain a closed manifold $\overline{P}$ as follows.   For each $j$, any gluing of solid tori along the boundary components of $P$ defines a Dehn filling of $H_j=\mathrm{cl}(P-N(\kappa_j))$ along all of its boundary components except $\partial N(\kappa_j)$.  Since each $H_j$ is hyperbolic, by Thurston's Hyperbolic Dehn Surgery Theorem \cite{BP,Th}, all but finitely many such Dehn fillings of $H_j$ result in a hyperbolic 3-manifold.  Furthermore, since $P$ is obtained by removing solid tori from $S^3$, for any integer $q$, if we attach the solid tori $Y_1$, \dots, $Y_{n+2}$ to $P$ with slope $\frac{1}{q}$, then $\overline{P}=S^3$. In this case each $H_j\cup Y_1\cup \dots\cup Y_{n+2}$ is the complement of a knot in $S^3$.  There are only finitely many $H_j$'s, and for each $j$, only finitely many slopes $\frac{1}{q}$ are excluded by Thurston's Hyperbolic Dehn Surgery Theorem.  Thus there is some integer $q$ such that if we glue the solid tori $Y_1$, \dots, $Y_{n+2}$ to any of the $H_j$ along $\partial N(C_1)$, \dots, $\partial N(C_{n})$, $\partial N(J)$, $\partial V$ with
slope  $\frac{1}{q}$, then we obtain the complement of a hyperbolic knot in $S^3$.  

Let $\Gamma_2$ denote the re-embedding of $\Gamma_1$, obtained as a result of gluing the solid tori $Y_1$, \dots, $Y_{n+2}$ to the boundary components of $P$ with slope $\frac{1}{q}$.  Now $\Gamma_2$ is an embedding of $K_n$ in $S^3$ such that every non-trivial knot in $\Gamma_2$ is hyperbolic.  Now there is a minor $G'$ of the embedded graph $\Gamma_2$ which is an embedding of our original graph $G$, such that every non-trivial knot in $G'$ is hyperbolic.
\end{proof}


\bigskip

\small

\begin{thebibliography}{10}

\bibitem{BP} R. Benedetti and C. Petronio, {\it Lectures on Hyperbolic Geometry}, Universitext, Springer-Verlag, Berlin (1992). 

\bibitem{CG} J. Conway and C. Gordon,  {\it Knots and links in spatial graphs},  J. of Graph Theory, {\bf  7}, (1983),  445--453.

\bibitem{fl} E. Flapan,  {\it Intrinsic knotting and linking of 
complete graphs}, Algebraic and Geometric Topology, {\bf 2} (2002) 
371--380.

\bibitem{fmn} E. Flapan, B. Mellor, R. Naimi, {\it Intrinsic Knotting and Linking are Arbitrarily Complex}, to appear Fundamentica Mathematicae, ArXiV math/0610501.

\bibitem{Th} W. Thurston, {\it Three-dimensional Geometry and Topology},Vol. 1, edited by Silvio Levy, Princeton Mathematical Series, {\bf 35}, Princeton University Press (1997).

\bibitem{Th2} W. Thurston, {\it Three dimensional manifolds, Kleinian groups, and Hyperbolic Geometry}, Bull. Amer. Soc. {\bf 6} (1982) 357--381.

\bibitem{Wa} F. Waldhausen, {\it Eine Klasse von 3-dimensionalen Mannigfaltigkeiten I, II}, Invent. Math. {\bf 3} 308--333, (1967), ibid {\bf 4} 87--117, (1967).

\end{thebibliography}
 \end{document}